\documentclass[10pt]{article}
\usepackage{amssymb}

\usepackage{amsmath,amssymb}\usepackage{cite}\usepackage{mathrsfs}\usepackage[amsmath,thmmarks]{ntheorem}
\usepackage{listings}
\usepackage[titletoc]{appendix}\allowdisplaybreaks

\usepackage{ulem}

\makeatletter
\renewcommand\theequation{\thesection.\@arabic\c@equation}
\makeatother

\newtheorem{thm}{ Theorem}[section]%
\newtheorem{lem}[thm]{ Lemma}%

\newtheorem{Remark}[thm]{ Remark}%
\newtheorem{Pro}[thm]{ Proposition}%
\newtheorem{df}[thm]{Definition}%
\input cyracc.def

\begin{document}
\title{The number of fuzzy subgroups of a finite abelian group of order $p^{n}q^{m}$}
\author{Lingling Han and Xiuyun Guo}

\date{}
\maketitle

\baselineskip=16pt

\vskip0.5cm

{\bf Abstract}. The purpose of this paper is to determine the number of fuzzy subgroups of a finite abelian group of order $p^{n}q^{m}$. As an application of our main result, explicit formulas for the number of fuzzy subgroups of $\mathbb{Z}_{p}^{n}\times\mathbb{Z}_{q}^{m}$ and $\mathbb{Z}_{p^{n}}\times\mathbb{Z}_{q}^{m}$ are given.


\section {Introduction}
\mbox{\hskip0.21in}All groups considered in this paper are finite.

One of the most important problems of fuzzy group theory is to classify the fuzzy subgroups of an abelian group. This topic has enjoyed a rapid evolution in the last years and many interesting and important results have been given. Since this problem is quite complex, one has to treat the case of cyclic groups first. For example, the authors in \cite{2-8} determine the number of distinct fuzzy subgroups of a cyclic group of square-free order, and the authors in \cite{2-9,2-10,2-11,2-17} deal with the number for cyclic groups of order $p^nq^m$($p, q$ primes). The big progress has been made by T\u{a}rn\u{a}uceanu and Bentea in \cite{2}, they give an explicit formula for the number of fuzzy subgroups of a cyclic group. In addition, they establish a recurrence relation verified by the number of fuzzy subgroups of an elementary abelian $p$-group, which is available to count the fuzzy subgroups of an elementary abelian $p$-group. After these remarkable works, we should restate the following as T\u{a}rn\u{a}uceanu and Bentea said in \cite{2}.

{\it Open problem:} Although the problem of counting the distinct fuzzy subgroups was solved for finite cyclic groups
and finite elementary abelian groups, it remains still open for arbitrary finite abelian groups (or for other more large
classes of finite groups). This may constitute the subject of some future studies.

It is because of the above problem that researchers begin to focus on more general abelian groups. For example, Ngcibi, Murali and Makamba in \cite{6-3}
obtain a formula for the number of fuzzy subgroups of the group $\mathbb{Z}_{p^{m}}\times \mathbb{Z}_{p^{n}}$ for $n=1,2,3$, which has been extended by Oh in \cite{6-1} for all values of $n$.
In the present paper we are concerned with the number of fuzzy subgroups of a new class of abelian groups, which is the direct product of a Sylow $p$-subgroup and a Sylow $q$-subgroup, and an explicit formula for this number for the above new class of abelian groups is obtained. The interesting is that we establish a one-to-one correspondence between
a family of some kind of subgroup chains for an abelian group and a set of some kind of vector pairs when we prove our main result, which make us count the cardinal number of the set of some kind of vector pairs instead of counting the cardinal number of the family of some kind of subgroup chains for an abelian group. As an application of our main result, we give two explicit formulas for the special cases.

\section {Preliminaries}

\mbox{\hskip0.21in}In this section we recall some basic notions and results which are useful in the sequel. We also make some agreements. The terminologies and notations not mentioned here agree with standard usage(for example, see\cite{001,01,02,2-5,2-6}).

\begin{df}\label{d-1}{\rm\cite{2-5}}
~Let $S$ be a set, then a mapping $\mu:S\rightarrow[0,1]$ is called a fuzzy subset of $S$.
\end{df}

\begin{df}\label{d-3}{\rm\cite{2-5}}
~Let $S$ be a set and $\mu:S\rightarrow[0,1]$ be a fuzzy subset of $S$. For each $\alpha\in[0,1]$, we define the level subset of $\mu$ as follows:
$$_{\mu}S_{\alpha}=\{x\in S\mid\mu(x)\geq\alpha\}.$$
\end{df}

\begin{df}\label{d-2}{\rm\cite{2-5}}
~Let $G$ be a group and $\mu:G\rightarrow[0,1]$ be a fuzzy subset of $G$. We say that $\mu$ is a
fuzzy subgroup of $G$ if it satisfies the next two conditions:

$(1)$~$\mu(xy)\geq min\{\mu(x),\mu(y)\}$, for all $x,y\in G;$

$(2)$~$\mu(x^{-1})\geq\mu(x)$, for any $x\in G.$
\end{df}

Note that $\mu(x^{-1})=\mu(x)$, for any $x\in G$, and $\mu(e)=max~\mu(G)$, where $e$ is the identity of $G$. The following important theorem was proved in \cite{1-1}. It could help to characterize the fuzzy subgroup of a group $G$ by using the level subsets.
\begin{thm}\label{t-1}{\rm\cite[Theorem 2.1 and Theorem 2.2]{1-1}}
~Let $G$ be a group and $\mu:G\rightarrow[0,1]$ be a fuzzy subset of $G$. Then $\mu$ is a fuzzy subgroup of $G$ iff all the level subsets of $\mu$ are subgroups of $G$. 
\end{thm}

The fuzzy subgroups of a group can be classified up to some natural equivalence relations. One of widely used equivalence relations(used in \cite{2-7,2-8,2-9,2-10,2-11,6-3,6-1}) is shown below: Two fuzzy subgroups $\mu$ and $\eta$ of a group $G$ are equivalent and denoted by $\mu\approx\eta$(written as $\mu\sim\eta$ in some papers), if and only if they satisfy $\mu(x)>\mu(y)\Longleftrightarrow\eta(x)>\eta(y)~\text{for all}~x,y\in G$~and~$\mu(x)=0\Longleftrightarrow\eta(x)=0~\text{for all}~x\in G.$

Another equivalence relation $\sim$ is used in \cite{2-17,2}: Two fuzzy subgroups $\mu$ and $\eta$ of a group $G$ are equivalent, written as $\mu\sim\eta$, if $\mu(x)>\mu(y)\Longleftrightarrow\eta(x)>\eta(y)~\text{for all}~x,y\in G.$  And two fuzzy subgroups $\mu,\eta$ of $G$ are called distinct if $\mu\nsim\eta$.

Notice that the relation $\sim$ is just the more general form of the relation $\approx$. In fact, the following result shows that one of these two numbers gives us another.

\begin{Pro}\label{t-2}{\rm\cite[Proposition 2.2]{6-1}}
~Let $G$ be a group. Then the number $n_{1}$ of fuzzy subgroups of $G$ under the relation $\sim$ and the number $n_{2}$ of fuzzy subgroups of $G$ under the relation $\approx$ satisfy the following equality:
$$n_{2}=2n_{1}-1.$$
\end{Pro}

Note that we prefer to use the equivalence relation $\sim$ to count the fuzzy subgroups of a group in this paper because it has a closely connection with the level subset. Now let $\mu$ and $\eta$ be two fuzzy subgroups of a group $G$. We take $\mu(G)=\{\alpha_{1},\alpha_{2},\cdots,\alpha_{n}\}$ and assume that $\alpha_{1}>\alpha_{2}>\cdots>\alpha_{n}.$ Then $\mu$ determines the following proper subgroup chain of $G$ that ends in $G$:
$$\Delta:~_{\mu}G_{\alpha_{1}}<~_{\mu}G_{\alpha_{2}}<\cdots<~_{\mu}G_{\alpha_{n}}=G.$$
In \cite{2-17} a necessary and sufficient condition for $\mu$ and $\eta$ to be equivalent with respect to $\sim$ was given.

\begin{thm}\label{t-2}{\rm\cite[Theorem 6]{2-17}}
~Let $G$ be a group. Let $\mu$ and $\eta$ be two fuzzy subgroups of $G$. Then $\mu\sim\eta$ if and only if $\mu$ and $\eta$ have the same set of level subsets. 
\end{thm}

That is to say $\mu\sim\eta$ if and only if $\mu$ and $\eta$ determine the same subgroup chains of type $\Delta$. It means that there exists a one-to-one correspondence between the set of distinct fuzzy subgroups of $G$ and the set of proper subgroup chains of $G$ which end in $G$.

It is easy to see that the number of distinct fuzzy subgroups of the trivial group $\{e\}$ is $1$. Thus in the following, we mainly focus on non-trivial groups.

\begin{df}\label{d1}
~Let $n$ be a positive integer and let $G_1,G_2,\cdots,G_n$ be subgroups of a group $G$ with $G_j < G_{j+1}$ for $j=1,2,\cdots,n-1$, then the following subgroup chain of $G$
$$\Gamma: G_{1}<G_{2}<\cdots<G_{n}$$
is called a proper subgroup chain of $G$. In this case the integer $n-1$ is called the \textbf{length} of the proper subgroup chain $\Gamma$, and the subgroups $G_{1}$ and $G_{n}$ are called the initial term and the terminal term of $\Gamma$.
\end{df}

Note that if $\Gamma$ is a subgroup chain of a group $G$ such that its initial term and terminal term are $G$, then $\Gamma$ is a proper subgroup chain of group $G$ with length $0$.

\begin{df}\label{d2}
~Let $\Gamma_1$ and $\Gamma_2$ be two proper subgroup chains of a group $G$. If $\Gamma_1$ and $\Gamma_2$ contain same subgroups of $G$, then we say $\Gamma_1=\Gamma_2$. Otherwise, we say $\Gamma_1$ and $\Gamma_2$ are different.
\end{df}

\begin{df}\label{d3}
~Let $G$ be a group and let $S$ be the family of all proper subgroup chains of $G$. Now set

$N(G)=\{\Gamma\in S\mid \text{the terminal term of}~\Gamma~\text{is}~G\};$

$H(G)=\{\Gamma\in S\mid \text{the terminal term of}~\Gamma~\text{is}~G~\text{and the initial term of}~\Gamma~\text{is not}~\{e\}\};$

$H_{i}(G)=\{\Gamma\in H(G)\mid \text{the length of}~\Gamma~\text{is}~i-1,~where~i~\text{is a positive integer}\}.$\\
We use $n(G),~h(G)$ and $h_{i}(G)$ to denote the cardinal numbers of $N(G),~H(G)$ and $H_{i}(G)$, respectively.
\end{df}

\begin{Remark}\label{r1}
$(1)$ As shown above, it is obvious that the number of distinct fuzzy subgroups of $G$ is $n(G)$.

$(2)$ It is clear that $H_{i}(G)\subseteq H(G) \subseteq N(G)$ for any group $G$ and for every positive integer $i$. Therefore $h_{i}(G)\leq h(G) \leq n(G)$ for any group $G$.

$(3)$ $n(G)=1$ and $h(G)=0$ if and only if $G$ is a trivial group, and $h_{1}(G)=1$ when $G$ is any group with $G\neq \{e\}$.
\end{Remark}

\begin{lem}\label{l1}
~Let $G$ be a group with $G\neq \{e\}$. Then $n(G)=2h(G)$.
\end{lem}

\begin{pf}
For any subgroup chain $\Gamma$ in $H(G)$ with
$$\Gamma: G_{1}<G_{2}<\cdots<G_{n},$$
we may construct a new subgroup chain $\Gamma^*$ by using $\Gamma$
$$\Gamma^*: \{e\}<G_{1}<G_{2}<\cdots<G_{n}$$
and set the family of the subgroup chains of $G$
$$H^{*}(G)=\{ \Gamma^*  | \Gamma \in H(G) \}.$$
It is easy to see that $H(G)\cap H^{*}(G)=\emptyset$, $|H(G)|=|H^{*}(G)|$ and $N(G)=H(G)\cup H^{*}(G)$.
Thus $n(G)=2h(G)$.
\end{pf}$\hfill{}\Box$

The next lemma follows immediately from the above definitions and Lemma $\ref{l1}$.

\begin{lem}\label{l2}
~Let $G$ be a group with $G\neq \{e\}$. Then
$$n(G)=2h(G)=2\sum\limits_{i=1}h_{i}(G).$$
\end{lem}

\section {The number of subgroup chains of an abelian group of order $p^{n}q^{m}$}

\mbox{\hskip0.21in} Now let $p$ and $q$ be different primes, and let $G$ be an abelian group of order $p^{n}q^{m}$ with positive integers $n$ and $m$. Then, by Sylow theorem, there exists a Sylow $p$-subgroup $A$ and a Sylow $q$-subgroup $B$ such that $G=A\times B$. If $H$ is a subgroup of $G$, then,
by \cite[Chapter III. Theorem 2.3]{001}, there exists a Sylow $p$-subgroup $A_{1}$ and a Sylow $q$-subgroup $B_{1}$ of $H$ such that $H=A_{1}\times B_{1}$ with $A_{1}\leq A,B_{1}\leq B$. It is clear that this kind of decomposition of $H$ is uniquely determined by $H$. Thus, if
\begin{equation}\label{001}
\Gamma: G_{1}<G_{2}<\cdots<G_{s}=G
\end{equation}
is a subgroup chain in $H(G)$, then there exists an unique subgroup $A_i$ of $A$ and an unique subgroup $B_i$ of $B$ such that $G_i=A_i\times B_i$ for every $i=1,2,\cdots,s$. Hence we have subgroup chains $A_1\leq A_2\leq\cdots\leq A_s=A$ and $B_1\leq B_2\leq\cdots\leq B_s=B$. If we remove redundant terms and the identity subgroup $\{e\}$ in these two subgroup chains, then we have the proper subgroup chain in $H(A)$ and the proper subgroup chain in $H(B)$ as follows:
\begin{equation}\label{002}
\Gamma_1: A_{i_{1}}<A_{i_{2}}<\cdots<A_{i_{k}}=A,
\end{equation}
\begin{equation}\label{003}
\Gamma_2: B_{j_{1}}<B_{j_{2}}<\cdots<B_{j_{l}}=B.
\end{equation}
Notice that the subgroup chains $(\ref{002})$ and $(\ref{003})$ are uniquely determined by the subgroup chain $(\ref{001})$, we may give the following definition.
\begin{df}\label{d5}
~Let $G$ be an abelian group of order $p^{n}q^{m}$ with different primes $p,q$ and positive integers $n,m$, and let $A$ be the Sylow $p$-subgroup of $G$ and $B$ be the Sylow $q$-subgroup of $G$ with $G=A\times B$. We call the subgroup chains $(\ref{002})$ and $(\ref{003})$ the \textbf{restriction} of the subgroup chain $(\ref{001})$ concerning subgroups $A$ and $B$. We also say that the subgroup chain $(\ref{001})$ is restricted by the subgroup chains $(\ref{002})$ and $(\ref{003})$.
\end{df}

By the above discussion, we state the lemma as follows:
\begin{lem}\label{l3}
~Let $G$ be an abelian group of order $p^{n}q^{m}$ with different primes $p,q$ and positive integers $n,m$, and let $A$ be the Sylow $p$-subgroup of $G$ and $B$ be the Sylow $q$-subgroup of $G$ with $G=A\times B$. Then we have

$(1)$~The subgroup chains $(\ref{001})$, $(\ref{002})$ and $(\ref{003})$ belong to $H_{s}(G)$, $H_{k}(A)$ and $H_{l}(B)$. If we set $b=max\{k,l\}$, then $s\in\{b,b+1,\cdots,k+l\};$

$(2)$~For any subgroup chain $\Gamma\in H(G)$, there exist two unique subgroup chains $\Gamma_{1}\in H(A)$ and $\Gamma_{2}\in H(B)$ such that $\Gamma_{1}$ and $\Gamma_{2}$ are the restriction of $\Gamma$ concerning subgroups $A$ and $B$;

$(3)$~Conversely, for any given subgroup chains $(\ref{002})$ and $(\ref{003})$, then it is clear that the following two subgroup chains in $H(G)$ are restricted by the subgroup chains $(\ref{002})$ and $(\ref{003})$:
\begin{equation*}\label{0011}
A_{i_{1}}\times \{e\}<A_{i_{2}}\times \{e\}<\cdots<A_{i_{k}}\times \{e\}<A_{i_{k}}\times B_{j_{1}}<\cdots<A_{i_{k}}\times B_{j_{l}}=A\times B=G,
\end{equation*}
\begin{equation*}\label{0012}
\{e\}\times B_{j_{1}}<\{e\}\times B_{j_{2}}<\cdots<\{e\}\times B_{j_{l}}<A_{i_{1}}\times B_{j_{l}}<\cdots<A_{i_{k}}\times B_{j_{l}}=A\times B=G.
\end{equation*}
\end{lem}

In fact, for any two given subgroup chains, there may exist many subgroup chains in $H(G)$ are restricted by them. Next lemma tell us how many subgroup chains are restricted by two given subgroup chains.

\begin{lem}\label{l4}
~Let $G$ be an abelian group of order $p^{n}q^{m}$ with different primes $p,q$ and positive integers $n,m$, and let $A$ be the Sylow $p$-subgroup of $G$ and $B$ be the Sylow $q$-subgroup of $G$ with $G=A\times B$. If
\begin{equation*}\label{004}
\Gamma_{1}: A_{1}<A_{2}<\cdots<A_{i}=A
\end{equation*}
and
\begin{equation*}\label{005}
\Gamma_{2}: B_{1}<B_{2}<\cdots<B_{j}=B
\end{equation*}
are subgroup chains in $H_{i}(A)$ and $H_{j}(B)$($1\leq i\leq n, 1\leq j\leq m$), respectively. Without loss of generality, we may assume that $i\geq j$. Then there exist $h(i,j)$ subgroup chains in $H(G)$ such that they are all restricted by $\Gamma_{1}$ and $\Gamma_{2}$, where
$$h(i,j)=\sum\limits_{k=0}\limits^{j}\dbinom{i+k}{i}\dbinom{i}{j-k}.$$
\end{lem}

\begin{pf}
First we use $U_{i,j,s}$ to denote the set of all subgroup chains in $H_{s}(G)$ such that they are all restricted by $\Gamma_{1}$ and $\Gamma_{2}$. Then it is clear that $\Gamma\in U_{i,j,s}$ if and only if

$(1)$ $\Gamma\in H_{s}(G)$ with $s\in\{i,i+1,\cdots,i+j\}$ and

$(2)$ $\Gamma$ is restricted by $\Gamma_{1}$ and $\Gamma_{2}$.\\
Next we use $V_{i,j,s}$ to denote the set of all vector pairs $(\textbf{a},\textbf{b})$ with $\textbf{a}=(a_{1},a_{2},\cdots,a_{i})$ and $\textbf{b}=(b_{1},b_{2},\cdots,b_{j})$ satisfying the following conditions:

(i) $a_{1}<a_{2}<\cdots<a_{i}$, and $a_{1},a_{2},\cdots,a_{i}\in\{1,2,\cdots,s\}$;

(ii) $b_{1}<b_{2}<\cdots<b_{j}$, and $b_{1},b_{2},\cdots,b_{j}\in\{1,2,\cdots,s\}$;

(iii) $\{a_{1},a_{2},\cdots,a_{i}\}\cup\{b_{1},b_{2},\cdots,b_{j}\}=\{1,2,\cdots,s\}$.\\
It should be noted that we say two vector pairs $(\textbf{a},\textbf{b})$ and $(\textbf{a}',\textbf{b}')$ are same if and only if  $\textbf{a}=\textbf{a}'$ and $\textbf{b}=\textbf{b}'$.

Now we claim that there is a one-to-one correspondence between $U_{i,j,s}$ and $V_{i,j,s}$.

In fact, let
\begin{equation*}\label{006}
\Gamma: G_{1}<G_{2}<\cdots<G_{s}=G
\end{equation*}
be a subgroup chain in $U_{i,j,s}$. Then, by using above methods, there exist $A_{1t}$ and $B_{2t}$ such that $G_{t}=A_{1t}\times B_{2t}$ for every $t\in\{1,2,\cdots,s\}$. And therefore $\Gamma$ become the following subgroup chain:
\begin{equation*}\label{007}
\Gamma: A_{11}\times B_{21}<A_{12}\times B_{22}<\cdots<A_{1s}\times B_{2s}=A\times B=G.
\end{equation*}
Without loss of generality, we may assume that $s=i+k$ with $0\leq k\leq j$. For the following subgroup chains,
\begin{equation*}
\Gamma_{3}: A_{11}\leq A_{12}\leq\cdots\leq A_{1s}=A
\end{equation*}
and
\begin{equation*}
\Gamma_{4}: B_{21}\leq B_{22}\leq\cdots\leq B_{2s}=B,
\end{equation*}
we may delete the identity subgroup $\{e\}$ first, and then delete redundant terms. We should notice that we always delete the right one for repeated two subgroups. For example, if $A_{1t}=A_{1t+1}$, we always delete $A_{1t+1}$. By using this agreement, we may have two subseries $x_{1},x_{2},\cdots,x_{i}$ and $y_{1},y_{2},\cdots,y_{j}$ of the series $1,2,\cdots,s$ such that the subgroup chain $\Gamma$ are restricted by the following two subgroup chains.
\begin{equation*}\label{008}
\Gamma_{5}: A_{1x_{1}}<A_{1x_{2}}<\cdots<A_{1x_{i}}=A,
\end{equation*}
\begin{equation*}\label{009}
\Gamma_{6}: B_{2y_{1}}<B_{2y_{2}}<\cdots<B_{2y_{j}}=B.
\end{equation*}
By Lemma $\ref{l3}$, we see that $A_{h}=A_{1x_{h}}$ for any $h\in\{1,2,\cdots,i\}$ and $B_{f}=B_{2y_{f}}$ for any $f\in\{1,2,\cdots,j\}$ and therefore
$$\Gamma_{1}=\Gamma_{5},~~~~\Gamma_{2}=\Gamma_{6}.$$
It is clear that
$$x_{1}<x_{2}<\cdots<x_{i}~~and~~\{x_{1},x_{2},\cdots,x_{i}\}\subseteq\{1,2,\cdots,s\};$$
$$y_{1}<y_{2}<\cdots<y_{j}~~and~~\{y_{1},y_{2},\cdots,y_{j}\}\subseteq\{1,2,\cdots,s\}.$$
Let $\textbf{x}=(x_{1},x_{2},\cdots,x_{i})$, $\textbf{y}=(y_{1},y_{2},\cdots,y_{j})$.
In order to prove that $(\textbf{x},\textbf{y})$ belongs to $V_{i,j,s}$, we only need to prove that $(\textbf{x},\textbf{y})$ satisfies the condition (iii). Indeed, if $\{x_{1},x_{2},\cdots,x_{i}\}\cup\{y_{1},y_{2},\cdots,y_{j}\}\neq\{1,2,\cdots,s\}$, then there exists a positive integer $r\in\{1,2,\cdots,s\}$ such that $r\notin\{x_{1},x_{2},\cdots,x_{i}\}$ and $r\notin\{y_{1},y_{2},\cdots,y_{j}\}$. Notice that $x_{1},x_{2},\cdots,x_{i}$ and $y_{1},y_{2},\cdots,y_{j}$ are the subseries of the series $1,2,\cdots,s$, we see that we must delete the both subgroups $A_{1r}$ in $\Gamma_{3}$ and $B_{2r}$ in $\Gamma_{4}$ when we delate redundant terms in the above. That means $A_{1r-1}=A_{1r}$ and $B_{2r-1}=B_{2r}$, which contradicts to that the following subgroup chain
\begin{equation*}\label{007}
\Gamma: A_{11}\times B_{21}<A_{12}\times B_{22}<\cdots<A_{1r-1}\times B_{2r-1}=A_{1r}\times B_{2r}<\cdots<A_{1s}\times B_{2s}=A\times B=G
\end{equation*}
is a proper subgroup chain.
Thus $\{x_{1},x_{2},\cdots,x_{i}\}\cup\{y_{1},y_{2},\cdots,y_{j}\}=\{1,2,\cdots,s\}$, and therefore $(\textbf{x},\textbf{y})\in V_{i,j,s}$. In this case it is easy to see that the following
$$\varphi: U_{i,j,s}\longrightarrow V_{i,j,s}$$
$$~~~~~~~~~~\Gamma\longmapsto (\textbf{x},\textbf{y})$$
is a map from  $U_{i,j,s}$ to $V_{i,j,s}$.

Furthermore, for any vector pair $(\textbf{a},\textbf{b})\in V_{i,j,s}$ with $\textbf{a}=(a_{1},a_{2},\cdots,a_{i})$ and $\textbf{b}=(b_{1},b_{2},\cdots,b_{j})$, we may find a subgroup chain $\Gamma_{\textbf{a},\textbf{b}}$ in $U_{i,j,s}$ such that $\varphi(\Gamma_{\textbf{a},\textbf{b}})=(\textbf{a},\textbf{b})$. In fact, we may construct the subgroup chain $\Gamma_{\textbf{a},\textbf{b}}$ of $G$ by using the following natural process. Since $\{a_{1},a_{2},\cdots,a_{i}\}\cup\{b_{1},b_{2},\cdots,b_{j}\}=\{1,2,\cdots,s\}$, it must be $1\in\{a_{1},a_{2},\cdots,a_{i}\}$ or $1\in\{b_{1},b_{2},\cdots,b_{j}\}$. So we may set $D_1$ by the following way.
\[D_{1}=\begin{cases}
A_{1}\times B_{1},&\text{if $a_{1}=1$ and $b_{1}=1$};\\
A_{1}\times\{e\},&\text{if $a_{1}=1$ and $b_{1}\neq1$};\\
\{e\}\times B_{1},&\text{if $a_{1}\neq1$ and $b_{1}=1$}.
\end{cases}\]
If $s=1$, then $i=j=1$ and therefore $D_{1}=A_{1}\times B_{1}=G$ is the subgroup chain $\Gamma_{\textbf{a}, \textbf{b}}$ we need. Now we assume that $D_{t-1}(2\leq t\leq s)$ is given and $D_{t-1}=X \times Y$, where $X\in\{ \{e\}, A_{1}, A_{2}, \cdots, A_{i}\}$ and $Y\in\{ \{e\}, B_{1}, B_{2}, \cdots, B_{j}\}$. Then we consider the case $t$. Also since $t\in\{a_{1},a_{2},\cdots,a_{i}\}$ or $t\in\{b_{1},b_{2},\cdots,b_{j}\}$, there exists $h\in\{1,2,\cdots,i\}$ such that $t=a_{h}$ or exists $f\in\{1,2,\cdots,j\}$ such that $t=b_{f}$. So we may set $D_t$ as the follows.
\[D_{t}=\begin{cases}
A_{h}\times B_{f},&\text{if $t=a_{h}=b_{f}$};\\
A_{h}\times Y,&\text{if $t=a_{h}\neq b_{f}$};\\
X\times B_{f},&\text{if $t=b_{f}\neq a_{h}$}.
\end{cases}\]
By using the above methods, we may have subgroups $D_{1},D_{2},\cdots,D_{s}$ of $G$, and therefore we have the following subgroup chain
$$\Gamma_{\textbf{a},\textbf{b}}: D_{1}<D_{2}<\cdots<D_{s}.$$
And it is easy to see that $D_{1}\neq\{e\},D_{s}=G$ and $D_{t}<D_{t+1}$ for any $t=1,2,\cdots,s-1$. Therefore we have $\Gamma_{\textbf{a},\textbf{b}}\in U_{i,j,s}$ and $\varphi(\Gamma_{\textbf{a},\textbf{b}})=(\textbf{a},\textbf{b})$. Now we may define the following map
$$\psi: V_{i,j,s}\longrightarrow U_{i,j,s}$$
$$~~~(\textbf{a},\textbf{b})\longmapsto \Gamma_{\textbf{a},\textbf{b}}$$
It is clear that $\varphi\psi((\textbf{a},\textbf{b}))=(\textbf{a},\textbf{b})$ and $\psi\varphi(\Gamma)=\Gamma$. So the claim is proved.

By the above claim, we have
$$h(i,j)=\sum\limits_{s=i}\limits^{i+j}|U_{i,j,s}|=\sum\limits_{s=i}\limits^{i+j}|V_{i,j,s}|.$$
Since we assume that $s=i+k$. Then
$$h(i,j)=\sum\limits_{k=0}\limits^{j}|V_{i,j,i+k}|.$$
Now we calculate the $|V_{i,j,i+k}|$. For any vector pair $(\textbf{a},\textbf{b})\in V_{i,j,s}$ with $\textbf{a}=(a_{1},a_{2},\cdots,a_{i})$ and $\textbf{b}=(b_{1},b_{2},\cdots,b_{j})$, it is easy to see that $(\textbf{a},\textbf{b})$ can be determined once $\{a_{1},a_{2},\cdots,a_{i}\}$ and $\{a_{1},a_{2},\cdots,a_{i}\}\cap\{b_{1},b_{2},$ $\cdots,b_{j}\}$ are given. Noticing that
$$|\{a_{1},a_{2},\cdots,a_{i}\}\cap\{b_{1},b_{2},\cdots,b_{j}\}|=i+j-s=i+j-(i+k)=j-k,$$
we see
$$|V_{i,j,i+k}|=\dbinom{i+k}{i}\dbinom{i}{j-k}.$$
Therefore
$$h(i,j)=\sum\limits_{k=0}\limits^{j}\dbinom{i+k}{i}\dbinom{i}{j-k}.$$
The lemma is proved.
\end{pf}$\hfill{}\Box$

\begin{lem}\label{l5}
~Let $G$ be an abelian group of order $p^{n}q^{m}$ with different primes $p,q$ and positive integers $n,m$, and let $A$ be the Sylow $p$-subgroup of $G$ and $B$ be the Sylow $q$-subgroup of $G$ with $G=A\times B$. Suppose that $\Gamma_{1}$ and $\Gamma_{2}$ are two subgroup chains in $H(A)$ and that $\Gamma_{3}$ and $\Gamma_{4}$ are two subgroup chains in $H(B)$. Now set

$R=\{\Gamma\in H(G)\mid \text{the restriction of}~\Gamma~\text{is}~\Gamma_{1}~\text{and}~\Gamma_{3}\},$

$W=\{\Gamma\in H(G)\mid \text{the restriction of}~\Gamma~\text{is}~\Gamma_{2}~\text{and}~\Gamma_{4}\}.$\\
If~$\Gamma_{1}$ is different from $\Gamma_{2}$ or $\Gamma_{3}$ is different from $\Gamma_{4}$, then $R\cap W=\emptyset$.
\end{lem}
\begin{pf}If not, then there is a subgroup chain denoted by $\Gamma_{*}$ such that $\Gamma_{*}\in R\cap W$. Thus $\Gamma_{*}$ is restricted by $\Gamma_{1}$ and $\Gamma_{3}$, meanwhile, is also restricted by $\Gamma_{2}$ and $\Gamma_{4}$. According to Lemma $\ref{l3}$, we have $\Gamma_{1}=\Gamma_{2}$ and $~\Gamma_{3}=\Gamma_{4}$, a contradiction. This finishes the proof of the lemma.
\end{pf}$\hfill{}\Box$


\begin{thm}\label{t+1}
~Let $G$ be an abelian group of order $p^{n}q^{m}$ with different primes $p,q$ and positive integers $n,m$, and let $A$ be the Sylow $p$-subgroup of $G$ and $B$ be the Sylow $q$-subgroup of $G$ with $G=A\times B$. Then
$$n(G)=2\left(\sum_{\substack{1\leq i\leq n\\ 1\leq j\leq m\\  i\geq j}}h(i,j)h_{i}(A)h_{j}(B)+\sum_{\substack{1\leq i\leq n\\ 1\leq j\leq m\\  i<j}}h(j,i)h_{i}(A)h_{j}(B)\right),$$
where
$$h(i,j)=\sum\limits_{k=0}\limits^{j}\dbinom{i+k}{i}\dbinom{i}{j-k},~~~\text{for}~i\geq j.$$
\end{thm}

\begin{pf}
According to Lemma $\ref{l2}$, we only need to count $h(G)$. And by Lemma $\ref{l3}$, we can see that
$$H(G)=\{\Gamma\in H(G)\mid\Gamma~\text{is restricted by subgroup chains}~\Gamma_{1}~\text{and}~\Gamma_{2},~\Gamma_{1}\in H(A),~\Gamma_{2}\in H(B)\}.$$
Then, according to Lemma $\ref{l4}$, Lemma $\ref{l5}$ and the symmetry of $H(A)$ and $H(B)$, we see
$$n(G)=2\left(\sum_{\substack{1\leq i\leq n\\ 1\leq j\leq m\\  i\geq j}}h(i,j)h_{i}(A)h_{j}(B)+\sum_{\substack{1\leq i\leq n\\ 1\leq j\leq m\\  i<j}}h(j,i)h_{i}(A)h_{j}(B)\right),$$where
$$h(i,j)=\sum\limits_{k=0}\limits^{j}\dbinom{i+k}{i}\dbinom{i}{j-k}~~~\text{for}~i\geq j.$$
The theorem is proved.
\end{pf}$\hfill{}\Box$

\section {The number of the distinct fuzzy subgroups of $\mathbb{Z}_{p}^{n}\times \mathbb{Z}_{q}^{m}$}

\mbox{\hskip0.21in}Let $p$ be a prime and $n$ a positive integer. We use $\mathbb{Z}_{p}^{n}$ to denote an elementary abelian $p$-group of order $p^{n}$. According to Theorem $\ref{t+1}$, in order to count the number of distinct fuzzy subgroups of the abelian group $\mathbb{Z}_{p}^{n}\times \mathbb{Z}_{q}^{m}$, what we only need to do is clarifying the number of subgroup chains of $\mathbb{Z}_{p}^{n}$. We recall a well-known result in group theory first.

\begin{lem}\label{l7}\rm\cite[Chapter III. Theorem 8.5]{001}
~Let $p$ be a prime and let $G$ be an elementary abelian $p$-group of order $p^{n}$ with $n\geq 1$. Then the number of subgroups of $G$ of order $p^{m}$$(1\leq m\leq n)$ is
$$\begin{bmatrix}
n \\
m
\end{bmatrix}_{p}=\frac{(p^{n}-1)(p^{n-1}-1)\cdots(p^{n-m+1}-1)}{(p^{m}-1)(p^{m-1}-1)\cdots(p-1)}.$$
\end{lem}

Lemma $\ref{l7}$ allows us to get the formula for calculating the number of subgroup chains of elementary abelian $p$-groups.

\begin{lem}\label{l10}
~Let $G=\mathbb{Z}_{p}^{n}$ with $n\geq 1$. Then

~~~~$h_{1}(G)=1;$

~~~~$h_{2}(G)=\sum\limits_{i=1}\limits^{n-1}\begin{bmatrix}
n \\
i
\end{bmatrix}_{p};$

~~~~$h_{3}(G)=\sum\limits_{1\leq i_{1}<i_{2}\leq n-1}\begin{bmatrix}
n \\
i_{2}
\end{bmatrix}_{p}\begin{bmatrix}
i_{2} \\
i_{1}
\end{bmatrix}_{p};$

~~~~~~$\vdots$

~~~~$h_{k}(G)=\sum\limits_{1\leq i_{1}<i_{2}<\cdots<i_{k-1}\leq n-1}\begin{bmatrix}
n \\
i_{k-1}
\end{bmatrix}_{p}\begin{bmatrix}
i_{k-1} \\
i_{k-2}
\end{bmatrix}_{p}\cdots\begin{bmatrix}
i_{2} \\
i_{1}
\end{bmatrix}_{p},~2\leq k\leq n.$
\end{lem}

\begin{pf}
It is clear that $h_{1}(G)=1.$ Let $\Gamma$ be a subgroup chain in $H_{k}(G)$ with $2\leq k\leq n$ as follows:
$$\Gamma:G_{1}<G_{2}<\cdots<G_{k}=G.$$
Then we can naturally have the following $k$-dimensional vector
$$\alpha=(|G_{1}|,|G_{2}|,\cdots,|G_{k}|).$$
For convenience, we call $\alpha$ the order vector of $\Gamma$. It is clear that
$$|G_{i}|~\big|~|G_{i+1}|~\text{for}~i=1,2,\cdots,k-1,~\text{and}~|G_{k}|=p^{n}.$$
Now set
$$\Lambda=\{\alpha\mid\alpha~\text{is an order vector of}~\Gamma, \Gamma\in H_{k}(G)\}$$
and for any $\alpha \in \Lambda$,
$$H_{k}^{\alpha}(G)=\{\Gamma\in H_{k}(G)\mid\text{the order vector of}~\Gamma~\text{is}~\alpha\}.$$
Then it is clear that
$$H_{k}(G)=\bigcup_{\alpha\in\Lambda}H_{k}^{\alpha}(G),$$
and
$$H_{k}^{\alpha}(G)\cap H_{k}^{\beta}(G)=\emptyset~\text{if}~\alpha\neq\beta.$$
Thus
$$h_{k}(G)=\sum_{\alpha\in\Lambda}|H_{k}^{\alpha}(G)|.$$
It is easy to see that
$$\Lambda=\{\alpha=(p^{i_{1}},p^{i_{2}},\cdots,p^{i_{k-1}},p^{n})\mid1\leq i_{1}<i_{2}<\cdots<i_{k-1}\leq n-1\}.$$
Noticing that the subgroups of an elementary abelian $p$-group are still elementary abelian $p$-groups, we see, by Lemma $\ref{l7}$, that
$$h_{k}(G)=\sum\limits_{1\leq i_{1}<i_{2}<\cdots<i_{k-1}\leq n-1}\begin{bmatrix}
n \\
i_{k-1}
\end{bmatrix}_{p}\begin{bmatrix}
i_{k-1} \\
i_{k-2}
\end{bmatrix}_{p}\cdots\begin{bmatrix}
i_{2} \\
i_{1}
\end{bmatrix}_{p},~2\leq k\leq n.$$
The lemma is proved.
\end{pf}$\hfill{}\Box$

According to Theorem $\ref{t+1}$ and Lemma $\ref{l10}$, it is easy to get the following result.

\begin{thm}\label{tA}
~Let $p$ and $q$ be different primes, and let $n$ and $m$ be positive integers. Also let $A=\mathbb{Z}_{p}^{n}$ and $B=\mathbb{Z}_{q}^{m}$ such that $G=A\times B$. Then
$$n(G)=2\left(\sum_{\substack{1\leq i\leq n\\ 1\leq j\leq m\\  i\geq j}}h(i,j)h_{i}(A)h_{j}(B)+\sum_{\substack{1\leq i\leq n\\1\leq j\leq m\\i<j}}h(j,i)h_{i}(A)h_{j}(B)\right),$$
where
$$h(i,j)=\sum\limits_{k=0}\limits^{j}\dbinom{i+k}{i}\dbinom{i}{j-k},~~~\text{for}~i\geq j;$$
$$h_{1}(A)=1,~~h_{i}(A)=\sum\limits_{1\leq l_{1}<l_{2}<\cdots<l_{i-1}\leq n-1}\begin{bmatrix}
n \\
l_{i-1}
\end{bmatrix}_{p}\begin{bmatrix}
l_{i-1} \\
l_{i-2}
\end{bmatrix}_{p}\cdots\begin{bmatrix}
l_{2} \\
l_{1}
\end{bmatrix}_{p},~2\leq i\leq n;$$

$$h_{1}(B)=1,~~h_{j}(B)=\sum\limits_{1\leq l_{1}<l_{2}<\cdots<l_{j-1}\leq m-1}\begin{bmatrix}
m \\
l_{j-1}
\end{bmatrix}_{q}\begin{bmatrix}
l_{j-1} \\
l_{j-2}
\end{bmatrix}_{q}\cdots\begin{bmatrix}
l_{2} \\
l_{1}
\end{bmatrix}_{q},~2\leq j\leq m.$$
\end{thm}





\section {The number of the distinct fuzzy subgroups of $\mathbb{Z}_{p^{n}}\times \mathbb{Z}_{q}^{m}$}

\mbox{\hskip0.21in}Let $p$ be a prime and $n$ a positive integer. We use $\mathbb{Z}_{p^{n}}$ to denote a cyclic group of order $p^{n}$. We recall the following well-known result in group theory.

\begin{lem}\label{l8}\rm\cite[Chapter I. Theorem 2.20]{001}
~Let $G=\langle g\rangle$ be a cyclic group of order $n$. Then there is a unique subgroup of $G$ of order $d$ for any positive divisor $d$ of $n$. And the subgroup of order $d$ is $\langle g^{k}\rangle$ with $k=\frac{n}{d}.$
\end{lem}

Using Lemma $\ref{l8}$, we can get the formula for calculating the number of subgroup chains of $\mathbb{Z}_{p^{n}}$.

\begin{lem}\label{l9}
~Let $p$ be a prime and $G=\mathbb{Z}_{p^{n}}$ with $n\geq 1$. Then

~~~~$h_{1}(G)=1;$

~~~~$h_{2}(G)=\dbinom{n-1}{1};$

~~~~$h_{3}(G)=\dbinom{n-1}{2};$

~~~~~~$\vdots$

~~~~$h_{k}(G)=\dbinom{n-1}{k-1},~1\leq k\leq n.$
\end{lem}

\begin{pf}
By Lemma $\ref{l8}$, we can see that all the subgroups of $G$ are $\{e\}, G_{1}, G_{2}, \cdots, G_{n}$, where $|G_{k}|=p^{k}$ for $k=1,2,\cdots,n$. And these subgroups satisfy
$$\{e\}<G_{1}<G_{2}<\cdots<G_{n}=G.$$
Now it easily is verified that

~~~~$H_{1}(G)=\{G_{n}\mid G_{n}~\text{is a subgroup chain}\};$

~~~~$H_{2}(G)=\{G_{i}<G_{n}\mid1\leq i\leq n-1\};$

~~~~~~$\vdots$

~~~~$H_{k}(G)=\{G_{i_{1}}<G_{i_{2}}<\cdots<G_{i_{k-1}}<G_{n}\mid 1\leq i_{1}<i_{2}<\cdots<i_{k-1}\leq n-1\}.$\\
Then we have

~~~~$h_{1}(G)=1;$

~~~~$h_{2}(G)=\dbinom{n-1}{1};$

~~~~$h_{3}(G)=\dbinom{n-1}{2};$

~~~~~~$\vdots$

~~~~$h_{k}(G)=\dbinom{n-1}{k-1},~1\leq k\leq n.$\\
The lemma is proved.
\end{pf}$\hfill{}\Box$

Now according to Theorem $\ref{t+1}$, Lemma $\ref{l10}$ and Lemma $\ref{l9}$, it is easy to get the formula for calculating the number of distinct fuzzy subgroups of the abelian group $\mathbb{Z}_{p^{n}}\times \mathbb{Z}_{q}^{m}$.

\begin{thm}\label{tB}
~Let $p$ and $q$ be different primes, and let $n$ and $m$ be positive integers. Also let $A=\mathbb{Z}_{p^{n}}$ and $B=\mathbb{Z}_{q}^{m}$ such that $G=A\times B$.
 Then
$$n(G)=2\left(\sum_{\substack{1\leq i\leq n\\ 1\leq j\leq m\\  i\geq j}}h(i,j)h_{i}(A)h_{j}(B)+\sum_{\substack{1\leq i\leq n\\1\leq j\leq m\\i<j}}h(j,i)h_{i}(A)h_{j}(B)\right),$$
where
$$h(i,j)=\sum\limits_{k=0}\limits^{j}\dbinom{i+k}{i}\dbinom{i}{j-k},~~~\text{for}~i\geq j;$$
$$h_{i}(A)=\dbinom{n-1}{i-1},~1\leq i\leq n;$$

$$h_{1}(B)=1,~~h_{j}(B)=\sum\limits_{1\leq l_{1}<l_{2}<\cdots<l_{j-1}\leq m-1}\begin{bmatrix}
m \\
l_{j-1}
\end{bmatrix}_{q}\begin{bmatrix}
l_{j-1} \\
l_{j-2}
\end{bmatrix}_{q}\cdots\begin{bmatrix}
l_{2} \\
l_{1}
\end{bmatrix}_{q},~2\leq j\leq m.$$
\end{thm}


\end{document}